\documentclass[12pt]{amsart}

\usepackage{amssymb}
\usepackage{amsmath}
\evensidemargin\oddsidemargin

\begin{document}
\setcounter{page}{1}

\newtheorem{PROP}{Proposition}
\newtheorem{REMS}{Remark}
\newtheorem{LEM}{Lemma}
\newtheorem{THE}{Theorem\!\!}
\newtheorem{COR}{Corollary \!\!}

\renewcommand{\theTHE}{}

\newcommand{\eqnsection}{
\renewcommand{\theequation}{\thesection.\arabic{equation}}
    \makeatletter
    \csname  @addtoreset\endcsname{equation}{section}
    \makeatother}
\eqnsection

\def\a{\alpha}
\def\b{\beta}
\def\B{{\bf B}} 
\def\FF{{\mathcal{F}}} 
\def\GG{{\mathcal{G}}}
\def\MM{{\mathcal{M}}}
\def\ca{c_{\a}}
\def\ka{\kappa_{\a}}
\def\coa{c_{\a, 0}}
\def\cua{c_{\a, u}}
\def\cL{{\mathcal{L}}} 
\def\Ea{E_\a}
\def\eps{{\varepsilon}} 
\def\esp{{\mathbb{E}}} 
\def\Ga{{\Gamma}} 
\def\G{{\bf \Gamma}} 
\def\e{{\rm e}}
\def\i{{\rm i}}
\def\K{{\bf K}}
\def\Ka{{\bf K}_\a}
\def\L{{\bf L}}
\def\lbd{\lambda}
\def\lcr{\left[}
\def\lpa{\left(}
\def\lva{\left|}
\def\pb{{\mathbb{P}}}
\def\rl{{\mathbb{R}}}
\def\prst{{\prec_{st}}}
\def\prcvx{{\prec_{cx}}}
\def\rpa{\right)}
\def\rcr{\right]}
\def\rva{\right|}
\def\W{{\bf W}}
\def\X{{\bf X}}
\def\XX{{\mathcal X}}
\def\U{{\bf U}_\a}
\def\V{{\bf V}_\a}
\def\Un{{\bf 1}}
\def\Z{{\bf Z}}
\def\A{{\bf A}}
\def\AA{{\mathcal A}}

\def\claw{\stackrel{d}{\longrightarrow}}
\def\elaw{\stackrel{d}{=}}
\def\qed{\hfill$\square$}

\title[Area of a spectrally positive stable processes]
      {The area of a spectrally positive stable process stopped at zero}

\author[Julien Letemplier]{Julien Letemplier}

\address{Laboratoire Paul Painlev\'e, Universit\'e Lille 1, Cit\'e Scientifique, F-59655 Villeneuve d'Ascq Cedex. {\em Email} : {\tt ju.letemplier@gmail.com}}

\author[Thomas Simon]{Thomas Simon}

\address{Laboratoire Paul Painlev\'e, Universit\'e Lille 1, Cit\'e Scientifique, F-59655 Villeneuve d'Ascq Cedex. Laboratoire de physique th\'eorique et mod\`eles statistiques, Universit\'e  Paris Sud, B\^atiment 100, F-91405 Orsay Cedex. {\em Email} : {\tt simon@math.univ-lille1.fr}}



\begin{abstract} An identity in law for the area of a spectrally positive L\'evy $\a-$stable process stopped at zero is established. Extending that of Lefebvre \cite{Le} for Brownian motion, it involves an inverse Beta random variable and the square of a positive stable random variable. This identity entails that the stopped area is distributed as the perpetuity of a spectrally negative L\'evy process, and is hence self-decomposable. We also derive a convergent series representation for the density, whose behaviour at zero is shown to be Fr\'echet-like. 
\end{abstract}

\maketitle

\section{Introduction and statement of the results}

Let $\{B_t, \, t\ge 0\}$ be a standard linear Brownian motion, starting from one, and let $T = \inf\{ t>0, \, B_t =0\}$ be its first hitting time of zero. The random variable
$$\AA\; =\; \int_0^T\! B_s\, ds$$
has been investigated by Lefebvre, who obtained in Theorem 2 of \cite{Le} the simple identity in law
\begin{equation}
\label{MainB}
\AA\; \elaw\; \frac{2}{9 \G_{1/3}}
\end{equation}
where, here and throughout, $\G_a$ stands for the Gamma random variable with density
$$\frac{x^{a-1}}{\Ga(a)}\, e^{-x}\, \Un_{\{x > 0\}}.$$ 
The identity (\ref{MainB}) is obtained as a consequence of the closed expression for the Laplace transform of the bivariate random variable $(T, \AA)$ in terms of the Airy function - see Theorem 1 in \cite{Le}. 
As observed in \cite{La} p. 402, this latter expression can be easily derived thanks to the Feynman-Kac formula. Notice that Airy functions also appear in the expression of the Laplace transform of many other Brownian areas - see \cite{J1}, whose laws are more complicated than (\ref{MainB}).

In this paper, our concern is to generalize (\ref{MainB}) to the random variables
$$\AA_\a\; =\; \int_0^T \!L^{(\a)}_s\, ds$$
where $\{L_t^{(\a)}, \, t\ge 0\}$ is a strictly $\a-$stable L\'evy process without negative jumps, starting from one, and $T = \inf\{ t>0, \, L_t^{(\a)} =0\}$ is its first hitting time of zero. 
Without loss of generality we choose the normalization
$$\esp[e^{-tL^{(\a)}_1}]= e^{t^{\alpha}},\qquad t \ge 0,$$
where $\alpha \in [1,2]$ is the self-similarity parameter. We refer e.g. to Chapter 3 in \cite{S} for more information on stable L\'evy processes and the above normalization. The boundary cases $\a =1,2$ correspond to the unit drift resp. the Brownian motion with variance $\sqrt{2},$ so that we have
\begin{equation}
\label{Maa}
\AA_1\;=\; \frac{1}{2}\quad\quad\mbox{and}\quad\quad \AA_2\; \elaw\; \frac{1}{9\G_{1/3}}\cdot
\end{equation}  
The above second identity, which is actually the precise statement of Theorem 2 in \cite{Le}, follows from (\ref{MainB}) and the self-similarity of Brownian motion. In the case $\a\in (1,2)$ the generator of the process $L^{(\a)}$ is non-local, and it seems unappropriate to appeal to Feynman-Kac formul\ae\, in order to obtain a tractable expression for the Laplace transform for $\AA_\a.$ Moreover, the absence of transition densities written in closed form prevents from using explicit computations as in the Gaussian case - see \cite{Le} and the references therein, to handle the random variable $\AA_\a.$ Instead, we will compute the fractional moments of $\AA_\a$ and exhibit a multiplicative identity in law. Introduce the Beta random variable $\B_{a,b}$ with density
$$\frac{\Ga(a+b)}{\Ga(a)\Ga(b)}\, x^{a-1}(1-x)^{b-1}\, \Un_{( 0,1)}(x),$$ 
and the positive $a-$stable random variable $\Z_a$ with Laplace transform 
$$\esp[e^{-\lbd \Z_a}]\; =\; e^{-\lbd^a}, \qquad \lbd \ge 0.$$
Our main observation is the following.
\begin{THE} With the above notation, one has the independent factorization
\begin{equation}
\label{MainL}
\AA_\a\; \elaw\; \lpa\frac{\a+1}{4}\rpa\,\times\, \Z_{\frac{2}{\a + 1}}^{\,2}\, \times\, \B_{\frac{1}{2},\frac{\a -1}{2(\a+1)}}^{-1}
\end{equation}
for every $\a\in(1,2).$
\end{THE}
Observe that (\ref{MainL}) is in accordance with the two boundary cases: when $\a =1$ the two random variables on the right-hand side boil down to one, whereas when $\a = 2$ the following identity obtained in Theorem 1 of \cite{TS2}: 
$$\Z_{\frac{2}{3}}^{\,2}\; \elaw\; \frac{4}{27} \; \G_{\frac{2}{3}}^{-1}\, \times\,\B_{\frac{1}{3}, \frac{1}{6}}^{-1},$$
combined with the elementary factorization $\G_a\elaw\G_{a+b}\times\B_{a,b},$ allows to recover the second identity in (\ref{Maa}). The proof of (\ref{MainL}) relies on an identification of the fractional moments of $\AA_\a$. The explicit computation of the latter in terms of Gamma functions - see (\ref{Mom}) below - is made possible by the strong Markov property and some exact results on the stable Kolmogorov process recently obtained in \cite{PS}.\\

The inverse Gamma random variable involved in Lefebvre's identity shares a number of distributional properties related to infinite divisibility. Recall that a non-negative random variable $X$ is said to be self-decomposable if the following identities hold
$$X\; \elaw\; cX\; +\; X_c$$
for every $c\in (0,1)$ with $X_c$ independent of $X$, or equivalently if its log-Laplace transform is of the type
$$-\log\esp[e^{-\lbd X}]\; =\; a\lbd \; +\; \int_0^\infty (1- e^{-\lbd x})\, \frac{k(x)}{x} \, dx, \qquad \lbd \ge 0$$
for some $a \ge 0$ and a non-negative, non-increasing function $k$ integrating $1\wedge x^{-1}.$ See again Chapter 3 in \cite{S} for details. The fact that all inverse Gamma random variables are self-decomposable can be observed either by a direct and non-trivial computation on the Laplace transform, or by Dufresne's celebrated identity for the perpetuity of a Brownian motion with drift - see Section 3.2 in \cite{BS1} for details and references. In this paper we will use the second approach and show the same property for $\AA_\a.$

\begin{COR} For every $\a\in(1,2),$ the random variable $\AA_\a$ is self-decomposable.
\end{COR}

More precisely, it will be proved that $\AA_\a$ is distributed as the perpetuity of a spectrally negative L\'evy process which drifts towards $+\infty.$ The latter background integrated L\'evy process turns out to be tightly connected with the dual process $-L^{(\a)}$ conditioned to stay positive, and also to the Fr\'echet distribution which is hidden in the factorization (\ref{MainL}) and can be viewed as another perpetuity - see Remark 3 below. \\

It is known that inverse Gamma distributions also satisfy a property which is more stringent than self-decomposability. The law of a non-negative random variable $X$ is called a generalized Gamma convolution ($X\in\GG$ for short) if it is the weak limit of an independent sum of Gamma random variables, or equivalently if its log-Laplace transform reads 
$$-\log\esp[e^{-\lbd X}]\; =\; a\lbd \; +\; \int_0^\infty \,\log \lpa \frac{x}{x+\lbd}\rpa U(dx), \qquad \lbd \ge 0,$$
for some $a \ge 0$ and some non-negative measure $U$ satisfying certain integrability conditions, which is called the Thorin measure of $X.$ We refer to Chapter 3 in \cite{Bd} for more details on this notion refining that of self-decomposability. The fact that inverse Gamma random variables are all in $\GG,$ with explicit Thorin measure, is a consequence of a computation on Bessel functions - see again Section 3.2 in \cite{BS1} for details and references. Although we cannot write the Laplace transform of $\AA_\a$ in any tractable way, we can prove the following result.

\begin{COR} For every $\a\in(1,2),$ the law of $\sqrt{\AA_\a}$ is a generalized Gamma convolution, with infinite Thorin measure.
\end{COR}

A conjecture by Bondesson - see \cite{Bd} p. 97 and also Section 7 in \cite{BJTP} - states that the $\GG-$property is stable by power transformation of order greater than one. If this conjecture is true, then Corollary 2 entails that $\AA_\a\in\GG$ as well, a reinforcement of Corollary 1. It will be shown in Section 3.1. below that $\AA_{5/3}\in\GG,$ but we do not know as yet how to handle the other values of $\a$. 
  
The closed expression for the fractional moments of $\AA_\a$ can be inverted in order to give a convergent series representation for its density $f_{\AA_\a}$. Throughout this paper, we will set $f_X$ for the density of an absolutely continuous random variable $X.$ 

\begin{COR} The density of $\AA_\a$ has a convergent series representation:
$$f_{\AA_\a}(x)\;= \;\Ga\lpa\frac{\a}{\a+1}\rpa\,\times\;\sum_{n=0}^\infty \frac{(-1)^n (\a+1)^{\frac{n+1}{\a +1}-1}x^{-\frac{n+1}{\a +1}-1}}{n! \, \Gamma(1- \frac{n+1}{\a +1})\Gamma(1- \frac{n+2}{\a +1})}, \qquad x >0.$$
\end{COR}

Observe that when $\a = 2$ the above summation is made over $n=3p$ only, and that further simplifications lead to
\begin{equation}
\label{Masu}
f_{\AA_2}(x)\;= \;\frac{x^{-4/3}}{3^{2/3}\Ga(1/3)}\; \sum_{p=0}^\infty\; \frac{(-1)^p (9x)^{-p}}{p!} \; =\; \frac{\Ga(2/3)\, x^{-4/3}\,e^{-\frac{1}{9x}}}{2\pi\, 3^{1/6}},
\end{equation}
which is the expression to be found in Theorem 2 of \cite{Le}. The above corollary also entails the first order asymptotics
$$f_{\AA_\a}(x)\;\sim\;\frac{(\a+1)^{\frac{1}{\a +1}-1}x^{-\frac{1}{\a +1}-1}}{\Gamma(\frac{\a-1}{\a +1})} \qquad\mbox{as $x\to+\infty,$}$$
which has, up to the multiplicative constant, the same speed as that of the density of the factor $\Z_{2/(\a + 1)}^{\,2}$ at infinity - see Formula (14.31) in \cite{S}. On the other hand, it does not seem possible to deduce from the above series representation the exact behaviour of $f_{\AA_\a}$ at zero. Nevertheless, using the identity (\ref{MainL}) we can show the following estimate.

\begin{COR} The asymptotic behaviour of the density of $\AA_\a$ when $x\to 0+$ is
\begin{equation}
\label{Moss}
f_{\AA_\a}^{}(x) \;\sim\; \ka\; x^{\frac{\alpha^2}{1-\alpha^2}}\,e^{-\ca \,x^{\frac{1}{1-\alpha}}},
\end{equation}
with
$$\ka\; =\; \frac{\Ga(\frac{\a}{\a+1})\, \sqrt{\frac{\a+1}{\a-1}}}{2\pi\,(\a+1)^{\frac{\a}{\a^2-1}}}\qquad\quad \mbox{and}\qquad\quad \ca\; =\; (\a-1)(\a+1)^{\frac{\a}{1-\a}}.$$
\end{COR}
This shows that the behaviour of $f_{\AA_\a}$ at zero is that of the generalized Fr\'echet density
$$f_{\ca^{\a-1}\G_{1/(\a+1)}^{1-\a}}(x)\; =\;{\tilde \kappa}_\a\; x^{\frac{\alpha^2}{1-\alpha^2}}\,e^{-\ca \,x^{\frac{1}{1-\alpha}}},$$ 
up to the normalizing constant
$${\tilde \kappa}_\a\; =\; \frac{(\a-1)^{\frac{-\a}{\a+1}}}{\Ga(\frac{1}{\a+1})\,(\a+1)^{\frac{\a}{\a^2-1}}}$$
which does not coincide with $\ka$ except for $\a =2.$ Observe also that making $\a=2$ on the right-hand side of (\ref{Moss}) yields the density in (\ref{Masu}). It should be possible to obtain a full asymptotic expansion of $f_{\AA_\a}$ at zero with our method - see Remark 5 below. But we have not adressed this issue, which is believed to be awfully technical, in the present paper.

\section{Proofs}

\subsection{Proof of the Theorem} To simplify the notation we will set $L = L^{(\a)}.$ Introducing the area process
$$A_t\; =\; \int_0^t L_s\, ds, \qquad t \ge 0,$$
recall that the bivariate process $X = \{(A_t, L_t), \, t\ge 0\}$ is strongly Markovian and denote by $\pb_{(x,y)}$ its law starting from $(x,y).$ Consider the stopping time
$$S\;=\;\inf\{t> 0,\; A_t=0\}$$
and observe that under $\pb_{(0,1)}$ one has a.s. $S > T, A_T > 0,$ and $L_S < 0.$
Setting $\{\FF_t, \, t \ge 0\}$ for the natural completed filtration of $X$ and applying the strong Markov property at $T$ entails that for every $s\in \rl,$ one has
\begin{eqnarray*}
\esp_{(0,1)}[ \vert L_S\vert^{s-1}]\;=\; \esp_{(0,1)}[\, \esp [ \vert L_S\vert^{s-1} \vert\FF_T]] &= & \esp_{(0,1)}[ \,\esp_{(A_T, 0)}[\vert L_S\vert^{s-1}]]\\
&=& \esp \,[\AA_\a^{\frac{s-1}{\a +1}}]\,\times\,  \esp_{(1,0)}[\vert L_S\vert ^{s-1}],
\end{eqnarray*}
possibly with infinite values on both sides, where the second equality follows from the absence of negative jumps for $L,$ and the third equality from the self-similarity of $L$ and $A$ with respective indices $1/\a$ and $1+1/\a.$ Applying Theorem B in \cite{PS} in the particular case $\rho = 1/\a$ (beware that we consider here the dual process, with no positive jumps), we get
$$\esp_{(0,1)}[ \vert L_S\vert^{s-1}]\;=\;\frac{\sin (\frac{\pi s}{\a +1})}{\sin (\frac{\pi \a s}{\a +1})}$$
and
$$\esp_{(1,0)}[\vert L_S\vert ^{s-1}]\; =\;  \frac{(\a +1)^{\frac{1-s}{\a +1}}\Ga(\frac{\a +2}{\a +1})\Ga (\frac{1-s}{\a +1})\sin (\frac{\pi}{\a+1})}{\Ga(\frac{s}{\a +1})\Ga (1-s)\sin (\frac{\pi \a s}{\a +1})}$$
for all $\vert s\vert < 1 + 1/\a.$ Dividing and simplifying with the help of the complement formula for the Gamma function, we deduce
\begin{equation}
\label{Mom}
\esp \,[\AA_\a^s]\;=\;(\a +1)^s\,\times\,\frac{\Gamma(\frac{\alpha}{\alpha+1})\Gamma(1-(\a+1)s)}{\Gamma(\frac{\alpha}{\alpha +1} -s)\Gamma(1-s)}
\end{equation}
for all $s < 1/(\a +1).$ Applying the Legendre-Gauss multiplication formula for the Gamma function entails
$$\esp\,[\AA_\a^s]\;=\;\lpa\frac{\a +1}{4}\rpa^s\!\times\,\frac{\Gamma(1-(\a+1)s)}{\Gamma(1-2s)}\, \times\, \frac{\Gamma(\frac{\alpha}{\alpha+1})\Gamma(\frac{1}{2}-s)}{\Ga(\frac{1}{2})\Gamma(\frac{\alpha}{\alpha +1} -s)},$$
and we can conclude by a fractional moment identification, recalling (see e.g. Formula (25.5) in \cite{S} for the second expression) that
$$\esp\,[\B_{a,b}^s]\; =\; \frac{\Gamma(a+s)\Ga(a+b)}{\Ga(a)\Gamma(a+b+s)}\qquad \mbox{and}\qquad \esp\,[\Z_a^s]\; =\; \frac{\Gamma(1-\frac{s}{a})}{\Gamma(1-s)}\cdot$$
\qed

\begin{REMS} {\em (a) It is well-known and easy to see - see e.g. Theorem 46.3 in \cite{S} - that under $\pb_{(0,1)},$ the random variable $T$ is distributed as $\Z_{1/\a}.$ The above theorem provides hence a connection between $\AA_\a = A_T$ and the random variable $\Z_{2/(\a +1)}.$ Notice that one can also derive from (\ref{Mom}) another factorization:
$$\Z_{\frac{1}{\a +1}}\; \elaw\; (\a+1)^{-1}\, \times\,\G_{\!\!\frac{\a}{\a +1}}^{-1}\,\times\,  \AA_\a.$$
However, it seems difficult with our method to obtain some valuable information on the Mellin transform of the bivariate random variable $(T, A_T).$\\

(b) With the notation of our above proof, it is possible to derive the law of $A_T$ under $\pb_{(x,y)}$ for any $x\in\rl$ and $y >0,$ by the self-similarity of $L^{(\a)}.$
One finds
$$\AA_\a\; \elaw\; x\; +\;\lpa\frac{(\a+1)\,y^{\a+1}}{4}\rpa\,\times\, \Z_{\frac{2}{\a + 1}}^{\,2}\, \times\, \B_{\frac{1}{2},\frac{\a -1}{2(\a+1)}}^{-1}.$$}
\end{REMS}

\subsection{Proof of Corollary 1} Let us first observe that $\AA_\a$ is infinitely divisible, by a simple pathwise argument not relying on (\ref{MainL}). Setting $T_x= \inf \{t>0,\, L_t=x\}$ for all $x >0$ and using the fact that $L$ has no negative jumps, it is easy to see from the Markov property that under $\pb_{(0,1)},$ for every $n\ge 2,$ there is an independent  decomposition
$$\AA_\a\; =\; X_1^{(n)}\; +\; \cdots\; +\; X_n^{(n)}$$
where
$$X_i^{(n)}\; \elaw\; A^{}_{T_{\frac{n-i}{n}}}\;\;\mbox{under $\pb_{(0,\frac{n+1-i}{n})}$}$$
for every $i=1,\ldots, n.$ Moreover, one has $T_{1-1/n}\to 0$ a.s. under $\pb_{(0,1)}$ as $n\to +\infty$ (by the well-known fact - see e.g. Theorem 47.1 in \cite{S} - that $L$ visits immediately the negative half-line when starting from 0), so that $X_1^{(n)}\to 0$ a.s. under $\pb_{(0,1)}$ when $n\to +\infty$ as well. Last, it is straightforward that 
$$\pb[ \vert X_i^{(n)}\vert > \eps]\; \le\; \pb[ \vert X_1^{(n)}\vert > \eps]$$
for every $\eps > 0$ and $i=1,\ldots, n.$ Putting everything together and applying Khintchine's theorem on triangular arrays - see e.g. Theorem 9.3 in \cite{S}, entails that $\AA_\a$ is infinitely divisible.

\qed

\begin{REMS} {\em The above argument does not make use of the self-similarity of $L,$ and hence applies to any spectrally positive L\'evy process which is not a subordinator.}
\end{REMS}

We now proceed to the proof of the self-decomposability of $\AA_\a.$ We will use the same argument as in \cite{BS1, BS2}, expressing $\AA_\a$ as the perpetuity of a certain spectrally negative L\'evy process. Setting
$$\Psi_\a(u)\;= \frac{u \,\esp[\AA_\a^{-(u+1)}]}{ \esp[\AA_\a^{-u}]}$$
for every $u>0,$ we first deduce from (\ref{MainL}) the formula $\Psi_\a(u) =\Phi_\a(s),$ with the notation $s = (\a +1) u$ and
$$\Phi_\a(s)\;=\; \frac{\Gamma(\alpha+ s)}{\Gamma(s)}\cdot$$
Applying the Lemma in \cite{BS1} with $t =1$ shows, after some simplifications, that
$$\Phi_\a(s)\;=\;\Gamma(\alpha)\,s\;+\;\int_{-\infty}^0(e^{sx}-1-sx) f_{\alpha}(x)\,dx$$
with 
$$f_\a(x)\; =\;\frac{e^{\alpha x}}{\Gamma(-\alpha) (1-e^{x})^{\alpha +1}}\cdot$$
Since $f_\a$ integrates $x^2\wedge 1$ on $(-\infty,0),$ this entails that $\Psi_\a$ is the Laplace exponent of a spectrally negative L\'evy process. Applying Bertoin-Yor's criterion for perpetuities, as stated in \cite{BS2} pp. 8-9, shows that
\begin{equation}
\label{SNLP}
\AA_\a\; \elaw\; \int_0^\infty e^{-Z^{(\a)}_t}\, dt
\end{equation}
where $\{Z^{(\a)}_t, \, t\ge 0\}$ is the spectrally negative L\'evy process with Laplace exponent
$$\esp[e^{\lbd Z^{(\a)}_1}]\; =\; e^{\Psi_\a(\lbd)}\; =\; e^{\Phi_\a((\a +1)\lbd)}, \qquad \lbd \ge 0.$$
It is then easy to see from the representation (\ref{SNLP}) and the spectral negativity of $Z^{(\a)}$ that $\AA_\a$ is self-decomposable - see the end of the proof of the Theorem in \cite{BS1}.

\qed

\begin{REMS}{\em (a) The above expression (\ref{SNLP}) extends to the boundary cases $\a=1,2.$ When $\a = 1,$ the L\'evy process $Z^{(1)}$ has Laplace exponent $2\lbd,$ so that (\ref{SNLP}) boils down to
$$\AA_1\; \elaw\; \int_0^\infty e^{-2t}\, dt\; =\; \frac{1}{2}\cdot$$ 
When $\a =2,$ the L\'evy process $Z^{(2)}$ has Laplace exponent $3\lbd + 9\lbd^2,$ and (\ref{SNLP}) reads
$$\AA_2\; \elaw\; \int_0^\infty e^{3\sqrt{2}B_t -3t}\, dt\; \elaw\; \frac{1}{18}\,\int_0^\infty e^{B_t - t/6}\, dt\; \elaw\; \frac{1}{9\G_{1/3}},$$ 
the third identity in law being a particular case of Dufresne's identity. \\
 
(b) It follows from Corollary 2 in \cite{CC} that the spectrally negative L\'evy process $\{\xi^{\uparrow,n}_t,\, t\ge 0\}$ appearing in the Lamperti transform of the dual process $-L^{(\a)}$ conditioned to stay positive, has log-Laplace exponent
$$-\log \esp[e^{\lbd \xi^{\uparrow,n}_1}]\; =\; \Phi_\a(\lbd), \qquad \lbd \ge 0,$$
with the above notation for $\Phi_\a.$ This can be shown from the second formula in Corollary 2 of \cite{CC} written in an appropriate way, bewaring the unusual notation (7) therein for the negativity parameter and correcting a misprint (the $+$ before $c_-$ should be a $-$) in Formula (17) therein. This entails the curious identity in law
\begin{equation}
\label{Lamf}
A_T\; \elaw\; \int_0^\infty e^{- (\a+1) \xi^{\uparrow,n}_t}\, dt.
\end{equation}
Recall that Brownian motion conditioned to stay positive is distributed as a three-dimensional Bessel process, whose Lamperti process $\xi^{\uparrow, n}$ is the drifted Brownian motion $\{B_t +t/2, \, t\ge 0\}$ - see \cite{CC} p. 969 and the references therein, so that by Dufresne's identity we obtain
$$\int_0^\infty e^{- 3 \xi^{\uparrow,n}_t}\, dt\; \elaw\; \frac{2}{9\G_{\frac{1}{3}}}\; \elaw\; \AA.$$
This can be viewed as a particular case of (\ref{Lamf}), with the proper normalization. It is quite interesting to compare (\ref{Lamf}) with the identity
\begin{equation}
\label{LAM}
T\; \elaw\; \int_0^\infty e^{- \a \xi^{\uparrow,n}_t}\, dt,
\end{equation}
which follows from the above Remark 1 and Theorem 7 of \cite{CKP}. We do not have any sensible interpretation for the structural relationship between (\ref{Lamf}) and (\ref{LAM}). Notice last that the perpetuities of certain L\'evy processes with positive mean and jumping density  
$$\frac{K\,e^{b x}}{(1-e^{x})^{\alpha +1}}\, \Un_{\{ x<0\}}$$
for some $K, b >0$ and $\a\in (1,2)$ have been studied in \cite{Patou}. Observe in particular that the factorization obtained in Theorem 4.6 (2) therein shares some similarities with (\ref{MainL}).\\

(c) Combining the Kanter factorization - see Corollary 4.1 in \cite{K} - and (\ref{MainL}) shows the identity  
$$\AA_\a\; \elaw\; \lpa\frac{\a+1}{4}\rpa\,\times\,\G_1^{1-\a}\, \times\, \lpa \B_{\frac{1}{2},\frac{\a -1}{2(\a+1)}}\times\, \K_{\frac{2}{\a + 1}}^{\a+1}\rpa^{-1},$$
where $\K_a$ is the so-called Kanter random variable of index $a\in (0,1)$ - see Section 3 in \cite{TS2} for more details about this random variable. Consider now the Fr\'echet random variable $\G_1^{1-\a}$ appearing in the above factorization of $\AA_\a.$ In \cite{BS1}, it was shown that $\G_1^{1-\a}$ is also distributed as the perpetuity of a spectrally negative L\'evy process, and the latter turns out to be quite close to the above process $Z^{(\a)}.$ More precisely, it follows from the Lemma in \cite{BS1} with $t=1,$ the proof of the Theorem in \cite{BS1} and a change of variable that
\begin{equation}
\label{Frech}
\G_1^{1-\a}\; \elaw\; \int_0^\infty e^{-{\tilde Z}^{(\a)}_t}\, dt
\end{equation}
where $\{{\tilde Z}^{(\a)}_t, \, t\ge 0\}$ is the spectrally negative L\'evy process having log-Laplace exponent
$$-\log \esp[e^{\lbd {\tilde Z}^{(\a)}_1}]\; =\; (\a-1)^{-1} \Phi_\a((\a-1)\lbd)), \qquad \lbd \ge 0.$$
We do not have any explanation for this proximity between the two Laplace exponents of $Z^{(\a)}$ and ${\tilde Z}^{(\a)}.$ Notice that the identity (\ref{Frech}) can also been derived, in a complicated way, from Theorem 4.1 in \cite{Patou}.}
\end{REMS}

\subsection{Proof of Corollary 2} From (\ref{MainL}) and a general result by Bondesson - see Theorem 1 in \cite{BJTP}, in order to show $\sqrt{\AA_\a}\in\GG$ it is enough to prove that 
$$\B_{\frac{1}{2},\frac{\a -1}{2(\a+1)}}^{-1/2}\,\in\,\GG\qquad\mbox{and}\qquad\Z_{\frac{2}{\a+1}}\,\in\,\GG.$$
The second fact is well-known - see Example 3.2.1 in \cite{Bd}, whereas the first one follows at once from Theorem 2 in \cite{BS2}. Finally, the fact that $\sqrt{\AA_\a}$ has infinite Thorin measure follows from Theorem 4.1.4 in \cite{Bd} and the above formula (\ref{Mom}), which entails that $\sqrt{\AA_\a}$ has negative moments of all orders.

\qed 

\begin{REMS} {\em A combination of Corollaries 1 and 2 shows the following identities in law between three random integrals:
$$\int_0^\infty e^{-Z^{(\a)}_t}dt\; \elaw\;\int_0^T L^{(\a)}_t dt \; \elaw\; \lpa\int_0^\infty a^{(\a)}_t d\Ga_t\rpa^2,$$
where $\{\Ga_t, \, t\ge 0\}$ is the Gamma subordinator and $\{a^{(\a)}_t, \, t\ge 0\}$ some deterministic function which is related to the Thorin measure of $\sqrt{\AA_\a}$ - see Proposition 1.1 in \cite{JRY}.}
\end{REMS}

\subsection{Proof of Corollary 3} We will reason along the same lines as in Proposition 2 in \cite{TS1}, and omit some details. Applying the Mellin inversion formula yields first
$$f_{\AA_\a}(x)\;=\;\frac{1}{2\pi x} \int_{\mathbb{R}}\MM_\a(s)\, x^{-\i s} \,ds,$$
with the notation 
$$ \MM_\a(s)\;= \; \esp\,[\AA_\a^{\i s}]\; =\; (\a +1)^{\i s}\,\times\,\frac{\Gamma(\frac{\alpha}{\alpha+1})\Gamma(1-(\a+1)\i s)}{\Gamma(\frac{\alpha}{\alpha +1} -\i s)\Gamma(1-\i s)}$$
for every $s\in\rl.$
Suppose first $x > 1.$ We evaluate the above integral with the help of the residue theorem applied to the contour joining $-R$ to $R$ along the real axis, and $R$ to $-R$ along the half-circle plotted in the lower half-plane. It is easy to see that the integral along this half-circle vanishes as $R\to +\infty$, so that it remains to consider the singularities inside the big contour, which are located at $t_n= -\i (n+1)/(\alpha+1), n \geq 0.$
Computing
$$\text{Res}_{t_n}(\MM_\a(s) x^{-\i s})\; =\;-\i\,\Ga\lpa\frac{\a}{\a+1}\rpa\,\times\, \frac{(-1)^n (\a+1)^{\frac{n+1}{\a +1}-1}x^{-\frac{n+1}{\a +1}}}{n! \, \Gamma(1- \frac{n+1}{\a +1})\Gamma(1- \frac{n+2}{\a +1})}$$
we deduce
$$f_{\AA_\a}(x)\;= \;\Ga\lpa\frac{\a}{\a+1}\rpa\,\times\;\sum_{n=0}^\infty\; \frac{(-1)^n (\a+1)^{\frac{n+1}{\a +1}-1}x^{-\frac{n+1}{\a +1}-1}}{n! \, \Gamma(1- \frac{n+1}{\a +1})\Gamma(1- \frac{n+2}{\a +1})}$$
for every $x > 1,$ and hence for every $x >0$ by analytic continuation (Stirling's formula shows indeed that the series on the right-hand side converges absolutely for every $x > 0$). 

\qed

\subsection{Proof of Corollary 4} We will work on the random variable $\XX_\a\,=\,\AA_\a^{\frac{1}{1-\a}},$ in order to simplify the notation. Changing the variable, the required estimate is tantamount to
\begin{equation}
\label{MainAs}
f_{\XX_\a}^{}(x) \;\sim\; \frac{\Ga(\frac{\a}{\a+1})\,\sqrt{\alpha^2-1}}{2\pi\, (\a+1)^{\frac{\a}{\a^2-1}}}\; x^{-\frac{\a}{\a+1}}\,e^{-\ca x},\qquad x\to 0\!+\!.
\end{equation}
Evaluating with (\ref{Mom}) the positive entire moments
$$\esp[\XX_\a^n]\;=\;(\a +1)^{\frac{n}{1-\a}}\,\times\,\frac{\Gamma(\frac{\alpha}{\alpha+1})}{\Gamma(\frac{\alpha}{\alpha +1}+\frac{n}{\alpha-1})}\,\times \,\frac{\Gamma(1+\frac{(\a+1)n}{\alpha-1})}{\Gamma(1+\frac{n}{\alpha-1})}$$
for every $n\ge 0,$ and applying Stirling's formula shows that
$$\frac{(\esp[\XX_\a^n])^{\frac{1}{n}}}{n}\; \to\; \frac{(\a+1)^{\frac{\a}{\a-1}}}{\e (\a-1)}\qquad \mbox{as $n\to \infty.$}$$
By a theorem of Davies-Kasahara (see Corollary 4.12.5 in \cite{BGT}, or Lemma 3.2 in \cite{CSY} for a more appropriate formulation), we deduce that 
$$x^{-1}\log \pb [\XX_\a > x]\; \to \; -\ca\qquad \mbox{as $x\to +\infty.$}$$
This yields the required asymptotic behaviour, at the logarithmic scale, for the survival function of $\XX_\a.$  Moreover, writing down via Fubini's theorem the moment generating function
$$\esp[e^{x\XX_\a}]\;=\;\sum_{n=0}^{\infty} a_n\, x^n, \qquad x > 0, $$
with
\begin{eqnarray*}
a_n & = & (\a +1)^{\frac{n}{1-\a}}\,\times\,\frac{\Gamma(\frac{\alpha}{\alpha+1})\Gamma(1+\frac{(\a+1)n}{\alpha-1})}{n!\,\Gamma(\frac{\alpha}{\alpha +1}+\frac{n}{\alpha-1})\Gamma(1+\frac{n}{\alpha-1})}\\
& \sim & 
\frac{\Ga(\frac{\a}{\a+1})\sqrt{\alpha^2-1}}{2\pi\, (\a-1)^{\frac{1}{\a+1}}} \;\ca^{-n}\, n^{\frac{-\a}{\a +1}}\;\;\;\mbox{as $n\to\infty,$}\\
\end{eqnarray*} 
and applying Karamata's theorem for power series - see Corollary 1.7.3 in \cite{BGT} - shows that  
\begin{equation}
\label{Kara}
\esp[e^{\ca x\XX_\a}]\;\sim\; \frac{\Ga(\frac{1}{\a+1})\,\Ga(\frac{\a}{\a+1})\,\sqrt{\alpha^2-1}}{2\pi\, (\a-1)^{\frac{1}{\a+1}} (1-x)^{\frac{1}{\alpha+1}}}\qquad \mbox{as $x \to 1\!-\!.$}
\end{equation}
At this stage, it is worth mentioning that (\ref{Kara}) can be obtained from (\ref{MainAs}) by integration. However, it does not seem that we can infer the reverse inclusion without any further assumption, such as the existence of a meromorphic extension in the neighbourhood of $\ca$ for the moment generating function - see Theorem 4 in \cite{FGD}, or a monotonicity condition on $f_{\XX_\a}$ at infinity - see Theorem 4.12.11 in \cite{BGT}, which we both could not prove {\em a priori}.\\

In order to show (\ref{MainAs}) rigorously and finish the proof, we will use the following power transformation of (\ref{MainL}):
$$\XX_\a\; \elaw\; \lpa\frac{4}{\a+1}\rpa^{\frac{1}{\a-1}}\, \times\, \B_{\frac{1}{2},\frac{\a -1}{2(\a+1)}}^{\frac{1}{\a-1}}\,\times\, \Z_{\frac{2}{\a + 1}}^{\frac{2}{1-\a}}.$$
The multiplicative convolution formula and a change of variable entails
\begin{equation}
\label{Mult}
f_{\XX_\a}^{}(x)\; =\; \int_0^\infty \lpa \frac{1}{1+y}\rpa f_{\,\U}\!\lpa\frac{1}{1+y}\rpa\; f_{\,\V} (x + xy)\, dy,
\end{equation}
where we have set
$$\U\; =\; \B_{\frac{1}{2},\frac{\a -1}{2(\a+1)}}^{\frac{1}{\a-1}}\qquad\mbox{and}\qquad \V\; =\; \lpa\frac{4}{\a+1}\rpa^{\frac{1}{\a-1}}\! \times\,\Z_{\frac{2}{\a + 1}}^{\frac{2}{1-\a}}.$$
On the one hand, we have
$$\lpa \frac{1}{1+y}\rpa f_{\,\U}\!\lpa\frac{1}{1+y}\rpa\; =\;\frac{(\a-1)^{\frac{\a-1}{2(\a+1)}}\,\Ga(\frac{\a}{\a+1})}{\sqrt{\pi}\,\Ga(\frac{\a-1}{2(\a+1)})}\; y^{\frac{\a-1}{2(\a+1)}-1}(1+O(y))\qquad \mbox{as $y\to 0\!+\!.$}$$
On the other hand, Formula (14.35) in \cite{S} entail after a change of variable and several simplifications 
$$f_{\,\V} (z)\; =\;\frac{\sqrt{\a-1}\,(\a+1)^{\frac{\a^2}{2(\a^2-1)} - \frac{1}{\a-1}}}{2\sqrt{\pi}}\; z^{-1/2}\; e^{-\ca  z} (1+ O(z^{-1/2}))\qquad \mbox{as $z\to +\infty.$}$$
Plugging these two first order expansions in the integral (\ref{Mult}), and making further simplifications, yields finally the required asymptotic behaviour (\ref{MainAs}).

\qed

\begin{REMS} {\em The two asymptotic expansions for the above $f_{\,\U}$ and $f_{\,\V}$ can be continued at every order - see again Formula (14.35) in \cite{S} for the second function. This could be used to obtain a refined expansion for $f_{\XX_\a}$ at infinity, or equivalently for $f_{\AA_\a}$ at zero. Displaying the full asymptotic expansion of $f_{\AA_\a}$ at zero seems however to be a very painful task.}
\end{REMS}

\section{Final remarks}

\subsection{On the $\GG-$property} As mentioned in the introduction, the fact that $\sqrt{\AA_\a}\in \GG$ gives some credit to the $\GG-$property for $\AA_\a$ itself. This refinement of Corollary 1 would also shed some analytic light on the self-decomposability of $\AA_\a.$ We can show this property for $\a =5/3,$ where a combination of (\ref{MainL}) and Theorem 1 in \cite{TS2} entails the identity in law
$$\AA_{5/3}\; \elaw\; \lpa\frac{3}{2^{13/3}}\rpa\,\times\, \G_{\frac{3}{4}}^{-\frac{2}{3}}\, \times\, \B_{\frac{1}{4},\frac{1}{12}}^{-\frac{2}{3}}\, \times\, \B_{\frac{1}{2},\frac{1}{6}}^{-\frac{2}{3}}\, \times\, \B_{\frac{1}{2},\frac{1}{8}}^{-1}.$$
Indeed, all random variable on the right-hand side are in $\GG$ by Theorems 2 and 4 in \cite{BS2}, and we can conclude by Theorem 1 in \cite{BJTP}. Notice that the same kind of argument would show the property for all $\a\in (1,2)$ if we could prove that $\Z_\b^2\in\GG$ for all $\b\in(2/3,1)$ and that by the stability of the $\GG-$property with respect to weak convergence - see Theorem 3.1.5 in \cite{Bd}, it is enough to consider $\b$ rational. Again, Theorem 1 in \cite{TS2} shows a certain factorisation of $\Z_\b^2$ with $\b$ rational into the product of one Fr\'echet and several Beta random variables at a unique negative power, which unfortunately becomes too small when $\a \neq 5/3,$ so that cannot use Theorem 5.2 in \cite{BS2} as above. Observe also from Section 3.1 in \cite{BS2} that small negative powers of Beta random variables may not be in $\GG$.
To show the plausible fact that $\Z_\b^2\in\GG$ for all $\b\in(0,1)$ rational, it could be useful to write down the density of the Beta-Gamma product given in Formula (2.3) of \cite{TS2} as a Meijer $G$-function - see Theorem 9 in \cite{ST}, although it does not seem easy at first sight to express the Laplace transform in a tractable way.

\subsection{On the bell-shape property} It follows from Corollaries 3 and 4 that the density function of $\AA_\a$ is real analytic on $(0, +\infty),$ and that all its derivatives vanish at zero and at infinity. Moreover, a consequence of Corollary 2, Wolfe-Yamazato's theorem - see e.g. Theorem 51.3 in \cite{S}, and the principle of isolated zeroes, is the strict unimodality of this density function, that is its first derivative vanishes only once on $(0,+\infty).$ By Rolle's theorem, we deduce that its $n-$th derivative vanishes at least $n$ times on $(0,+\infty),$ and one can ask whether it vanishes exactly $n$ times for all $n\ge 1.$ Such a property, which is called the bell-shape in the literature, has been conjectured in \cite{TS3} for all positive self-decomposable distributions having an infinite spectral function at zero - see Conjecture 1 therein. Observe that $\AA_\a$ has such a self-decomposable law by Corollary 2, Corollary 4 and a theorem by Zolotarev - see Remark 28.6 in \cite{S}. Observe also that drawing the density with the help of the series representation of Corollary 3 and some plotting software exhibits the visual bell-shape property for $f_{\AA_\a}$, whose second derivative does seem to vanish only twice on $(0,+\infty).$

\subsection{Stable processes with negative jumps} The present paper studies the law of the random variable
$$\int_0^T \!L^{(\a)}_s\, ds$$
where $\{L_t^{(\a)}, \, t\ge 0\}$ is a spectrally positive stable process. It is natural to ask if (\ref{MainL}) could be extended in the presence of negative jumps, with $T$ defined as the first passage time below zero. In this case, it is known that the stable process, when starting positive, crosses zero by a negative jump with an explicit expression for the law of the undershoot - see e.g. Remark 42.18 in \cite{S} and the references therein. In order to apply our method, it would be hence necessary to find a closed expression of $\esp_{(x,y)}[\vert L_S\vert^{s-1}]$ with $x >0$ and $y <0,$ with the notation of our above proof. When $y=0,$ Theorem B in \cite{PS} provides a formula for these fractional moments in terms of the Gamma function, in full generality on the stable process. However, it seems hard to get a tractable formula when $y <0.$

\end{document}